\documentclass[12pt,a4paper]{article}
\usepackage[utf8]{inputenc}
\usepackage[T1]{fontenc}
\usepackage[english]{babel}

\usepackage{graphicx} 
\usepackage{textcomp}
\usepackage{listings}
\usepackage[space=true]{accsupp}
\newcommand{\copyablespace}{\BeginAccSupp{method=hex,unicode,ActualText=00A0}\EndAccSupp{}}

\usepackage{amsmath,amssymb,mathtools}
\usepackage{amsthm}

\newtheorem{example}{Example}

\usepackage[colorlinks=true, linkcolor=blue, urlcolor=blue, citecolor=blue, anchorcolor=blue]{hyperref}
\usepackage[table,dvipsnames]{xcolor}
\usepackage{fancyhdr}
\usepackage{longtable}
\usepackage{multicol}
\usepackage{enumerate}
\usepackage{enumitem}
\setlist[itemize]{leftmargin=5.5mm}
\usepackage{multirow}
\usepackage[space=true]{accsupp} 
\usepackage{listings}
\usepackage{caption}
\usepackage{subcaption}

\usepackage{xurl}
\usepackage{tikz}
\usepackage{pgfplots}
\usepgfplotslibrary{colorbrewer}
\pgfplotsset{compat=newest, cycle list/Set1-8}

\usepackage[numbers,sort&compress]{natbib}

\begin{document}

\title{The paxotopy method for finding isolated roots of a system of nonlinear equations} 

\maketitle

\begin{center}
{S\'andor Boz\'oki} \\[1mm]
{Institute for Computer Science and Control (SZTAKI),} \\
{Corvinus University of Budapest} \\
Budapest, Hungary \\
\texttt{bozoki.sandor@sztaki.hu} \\
\end{center}

\begin{abstract}
We focus on systems of nonlinear equations, where some equations are non-polynomial.
Building on the remarkable developments in solving polynomial systems during the recent decades,
and on a few extensions to specific examples of the transcendental case, a general framework
for finding isolated roots is proposed. The name \textbf{paxotopy} is proposed after
(1) \textbf{p}olynomial \textbf{a}ppro\textbf{x}imation of the transcendental functions, then
(2) solving the approximating polynomial system by the hom\textbf{otopy} continuation method, and finally
(3) the roots of the polynomial systems are used as starting points of a Newton's iteration to arrive at the solutions of the original system. These steps combine the advantages of local and global methods.

\end{abstract}

\section*{}

Solving a system of nonlinear equations is a challenging problem, even if all equations are polynomial.
However, the methodologies of solving systems of polynomial equations those have been developed during the
recent decades, are general and robust enough to be extendable to the non-polynomial case. 

We propose the name \textbf{paxotopy} after
\begin{itemize}
    \item (1) \textbf{p}olynomial \textbf{a}ppro\textbf{x}imation of the transcendental functions, then
    \item (2) solving the approximating polynomial system by the hom\textbf{otopy} continuation method, and finally
    \item (3) the roots of the polynomial systems are used as starting points of a Newton's iteration to arrive at the solutions of the original system.
\end{itemize}
All the three steps include a family of possible options (the polynomial approximation can be based on, e.g.,
Taylor series or Chebyshev polynomials; homotopy and Newton are again wide collections of algorithms).

Ji, Wu, Feng, Li and Qin \cite[Example 4.4]{JiWuFengLiQin2016} and Boz\'{o}ki \cite{Bozoki2020}
 applied Taylor approximation in the first step and homotopy method (HOM4PS-2, HOM4PS-3, respectively)
 in the second one.
Boyd \cite[Figure 20.1]{Boyd2014} considered a system of two transcendental equations of two variables,
 and applied Chebyshev polynomials in the first step, while several approaches are listed for the second step.

Although polynomial systems can be solved by several other algorithms (e.g., Gr\"obner bases) besides homotopy,
the proposed name paxotopy is kept for the methods that use homotopy continuation in the second step.
The reason behind is that there is no specific need for exact (e.g., symbolic) solutions of an approximating
polynomial system, numerical solutions are sufficient.

Paxotopy can utilize the additional knowledge of roots' location if it exists. If we are interested the
(possibly large) neighborhood of a given point, that is however too far to be usable as a starting point
of a Newton iteration, then Taylor approximation of the transcendental functions is a reasonable option in step 1.
If the search space is given in the form of a box, then its intervals and Chebyshev approximation
of the transcendental functions on these intervals is recommended.
In case there is no information at all on the possible location(s) of the roots, a systematic
 divide and conquer strategy is applied.

\begin{example}

Consider Carnahan's system \cite{Carnahan1964},\cite[page 308]{CarnahanGourdinLutherWilkes1969}:

\begin{eqnarray*}
f(x,y) = \frac{1}{2}  \sin(x y)  - \frac{y}{4 \pi} - \frac{x}{2} = 0 \\
g(x,y) = \left( 1- \frac{1}{4\pi}\right) (e^{2x} -e)+ \frac{e y}{\pi} -2 e x = 0 \\
\end{eqnarray*}

\begin{figure}[ht]
\begin{center}
\includegraphics[scale=0.2]{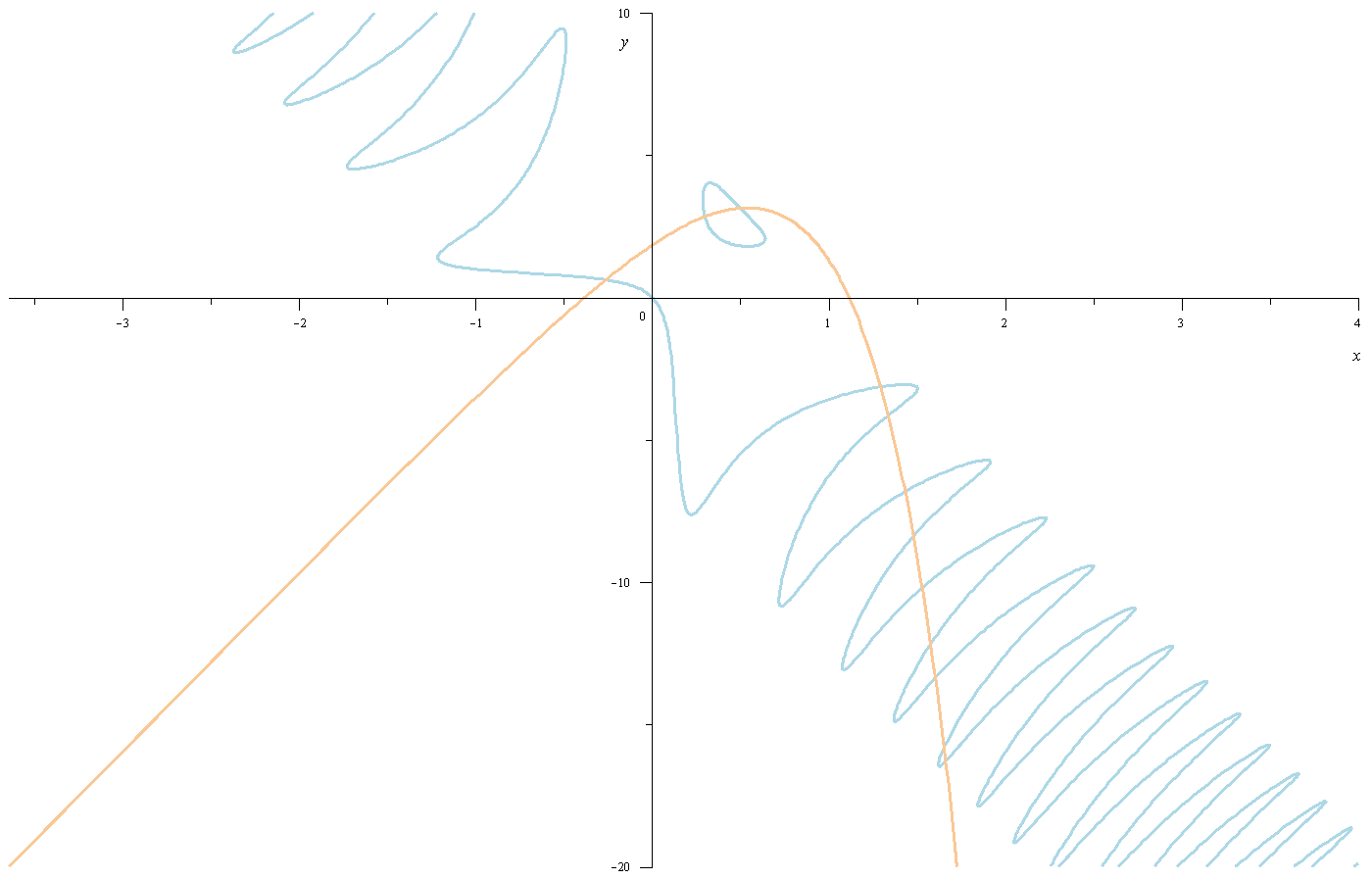}
\end{center}
\caption{Carnahan's system}
\end{figure}

Replace the $\sin$ function by its Chebyshev approximation in the interval $[-8,8]$ with the largest error $0.01$ and let the bivariate polynomial $p_1$ be 
as follows: \\

$p_1(x,y) = 
0.9881721447 xy - 0.1613036337 (xy)^3 + 0.007628424055 (xy)^5 \linebreak
- 0.0001576984074 (xy)^7 + 0.000001555955081 (xy)^9
 -6.10272221\cdot10^{-9} (xy)^{11} \linebreak - \frac{y}{4 \pi} - \frac{x}{2}. \\
$

Similarly, let $q_1(x,y)$ be the polynomial approximation of $g$ in $[-3,3].$
The solutions of $ p_1 = 0$ and $ q_1 = 0$ are plotted by dotted blue and red in Figure 2.
The green area $-8 \leq xy \leq 8, -3 \leq x \leq 3 $  shows the first search region.
Within this region, the solutions of $ f = 0$ and $ g = 0$ are well approximated, both pairs of curves overlap. Observe $p_1$ also has a branch in the negative orthant without having any part of $f$, but we will filter out such \emph{false} solutions.

\begin{figure}[ht]
\begin{center}
\includegraphics[scale=0.2]{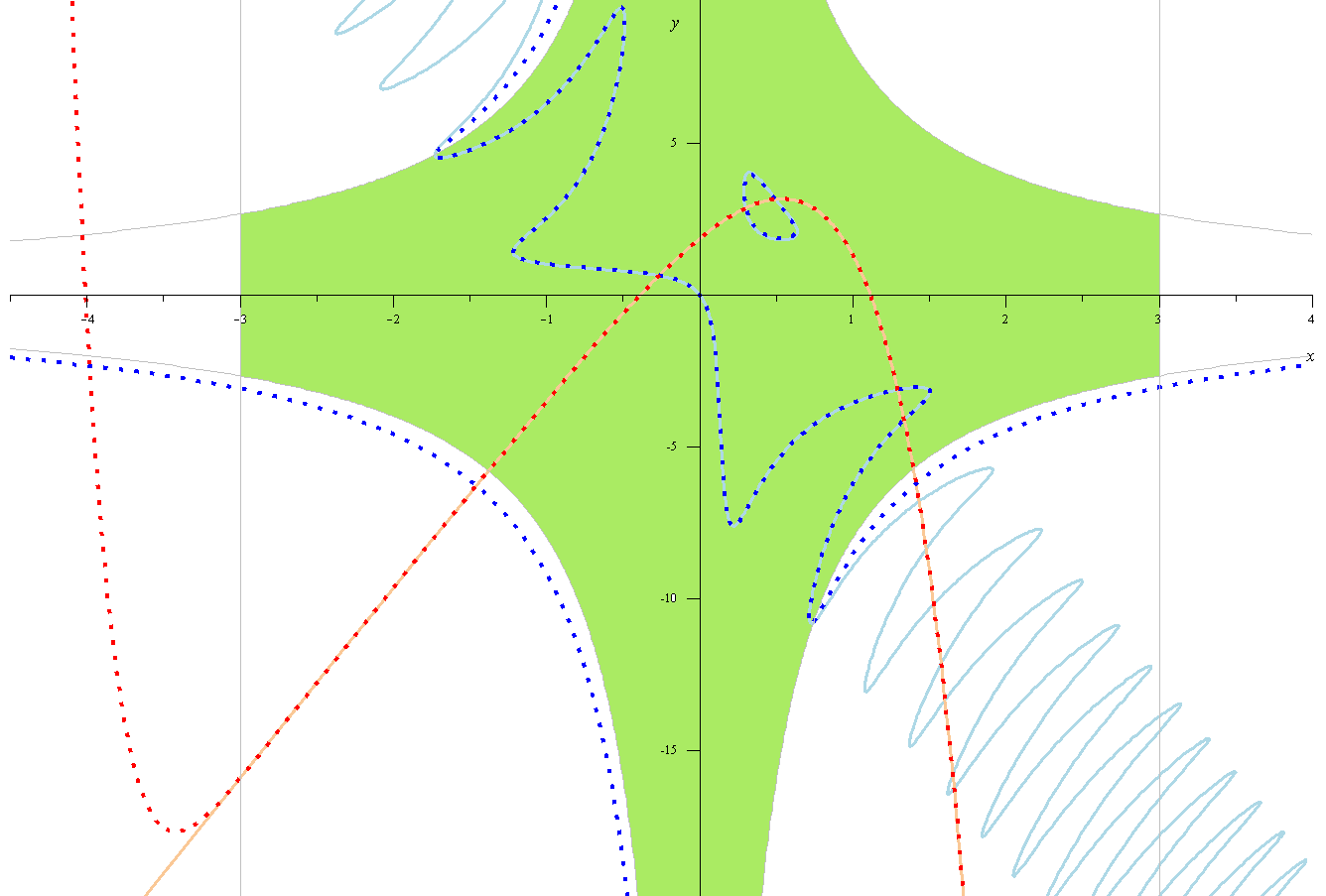}
\end{center}
\caption{Zeros of the approximating polynomials $p_1$ and $ q_1$}
\end{figure}

Figure 3 shows the solutions to the polynomial system $p_1(x,y) = 0, q_1(x,y) = 0$. Green dots indicate the ones with converging Newton's method written for the original system $f=g=0.$ Red dots correspond to solutions of the polynomial system only.

\begin{figure}[ht]
\begin{center}
\includegraphics[scale=0.2]{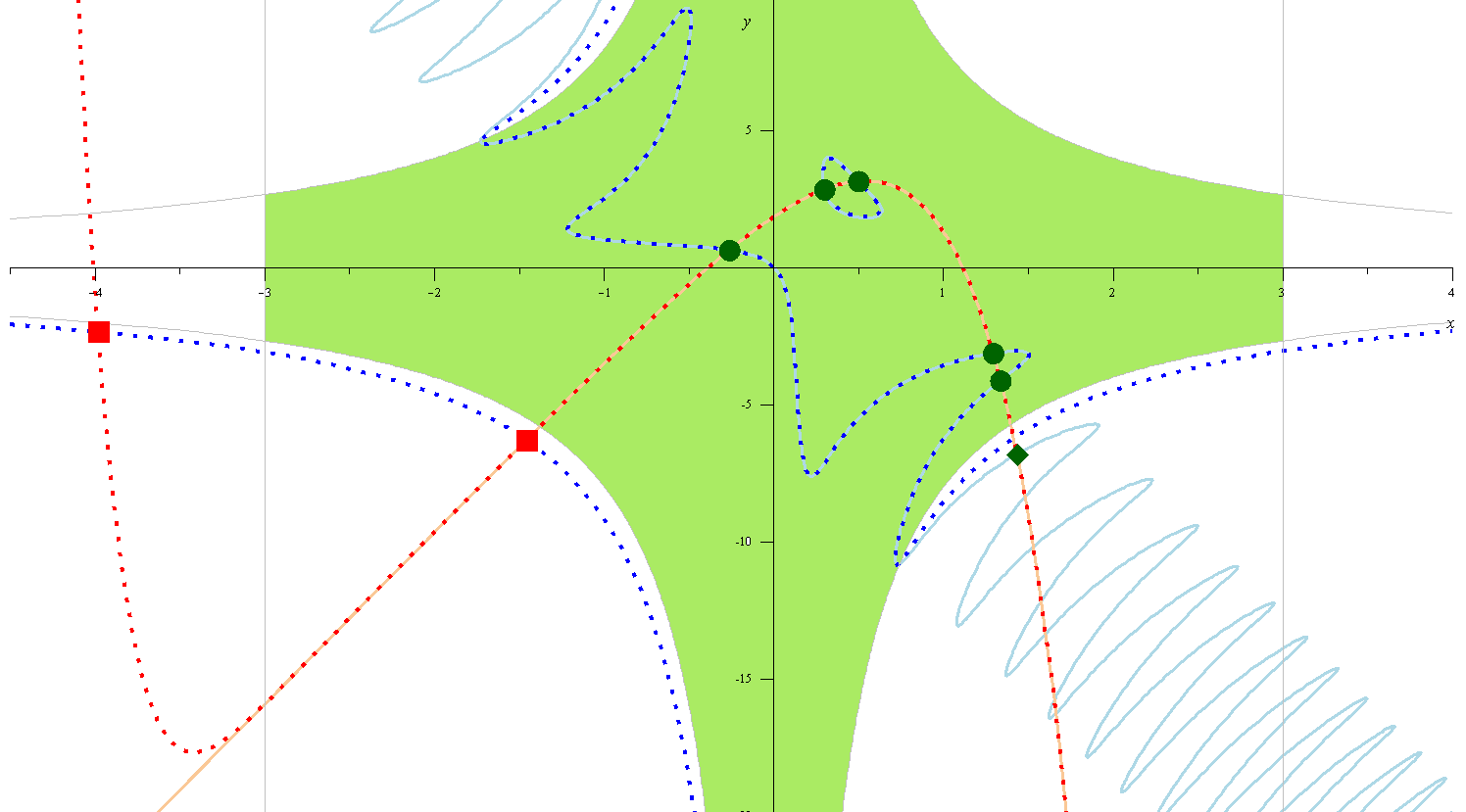}
\end{center}
\caption{Solutions to the polynomial system $p_1(x,y) = 0, q_1(x,y) = 0$}
\end{figure}

We can now choose the second search region as in Figure 4: 
take the Chebyshev approximation of $\sin$ in the interval $[8,24]$, 
leading to the polynomial approximation $p_2$ of $f$. 
Let $q_2$ be the same as $q_1$ (with the help of Figure 1). 

\begin{figure}[ht]
\begin{center}
\includegraphics[scale=0.2]{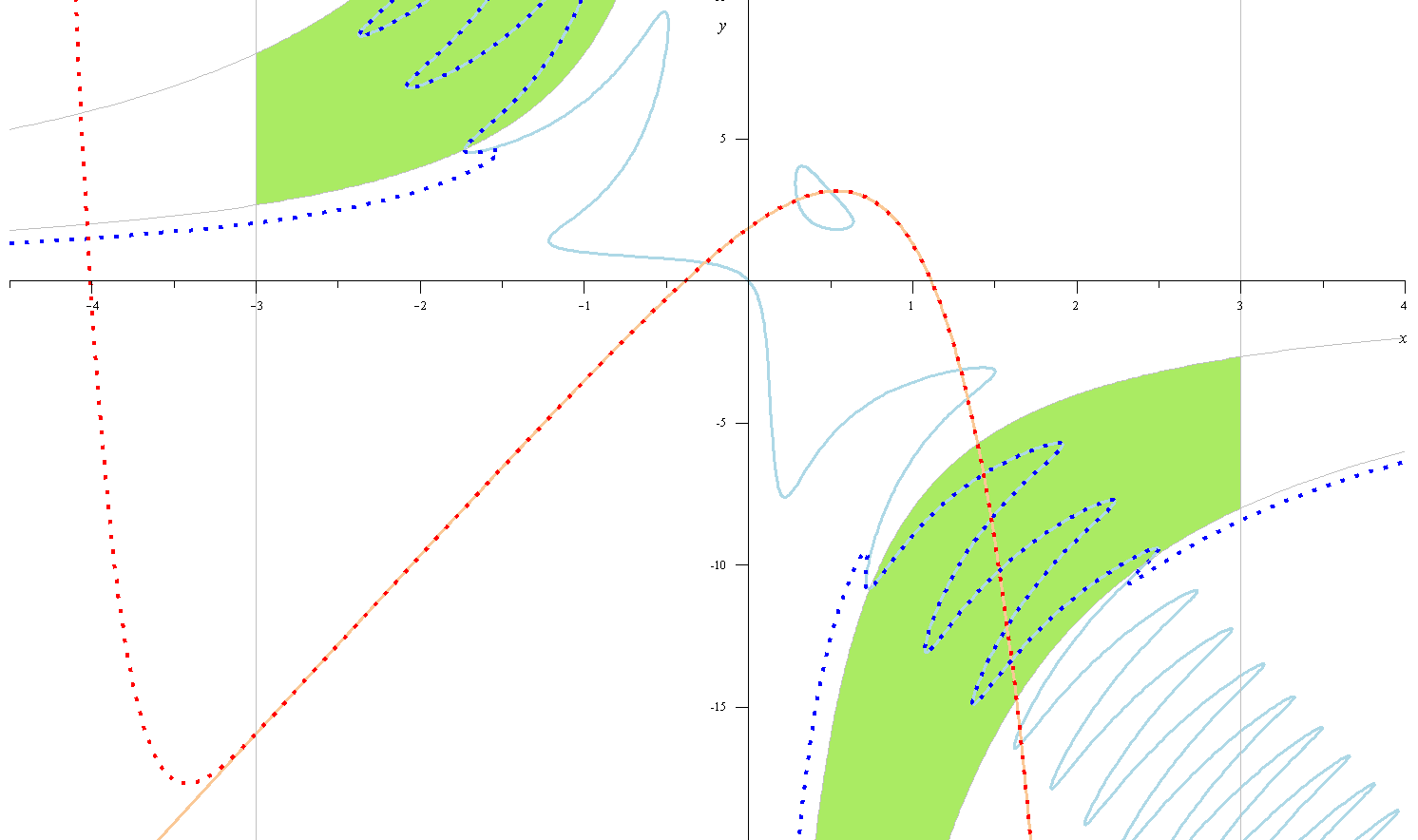}
\end{center}
\caption{Zeros of the approximating polynomials $p_2$ and $ q_2 (=q_1)$}
\end{figure}

Figure 5 shows the solutions to the polynomial system $p_2(x,y) = 0, q_2(x,y) = 0$. Green dots indicate the ones with converging Newton's method written for the original system $f=g=0.$ The red dot corresponds to a false solution again, that solves the polynomial system only, but not the original $f=g=0.$

\begin{figure}[ht]
\begin{center}
\includegraphics[scale=0.2]{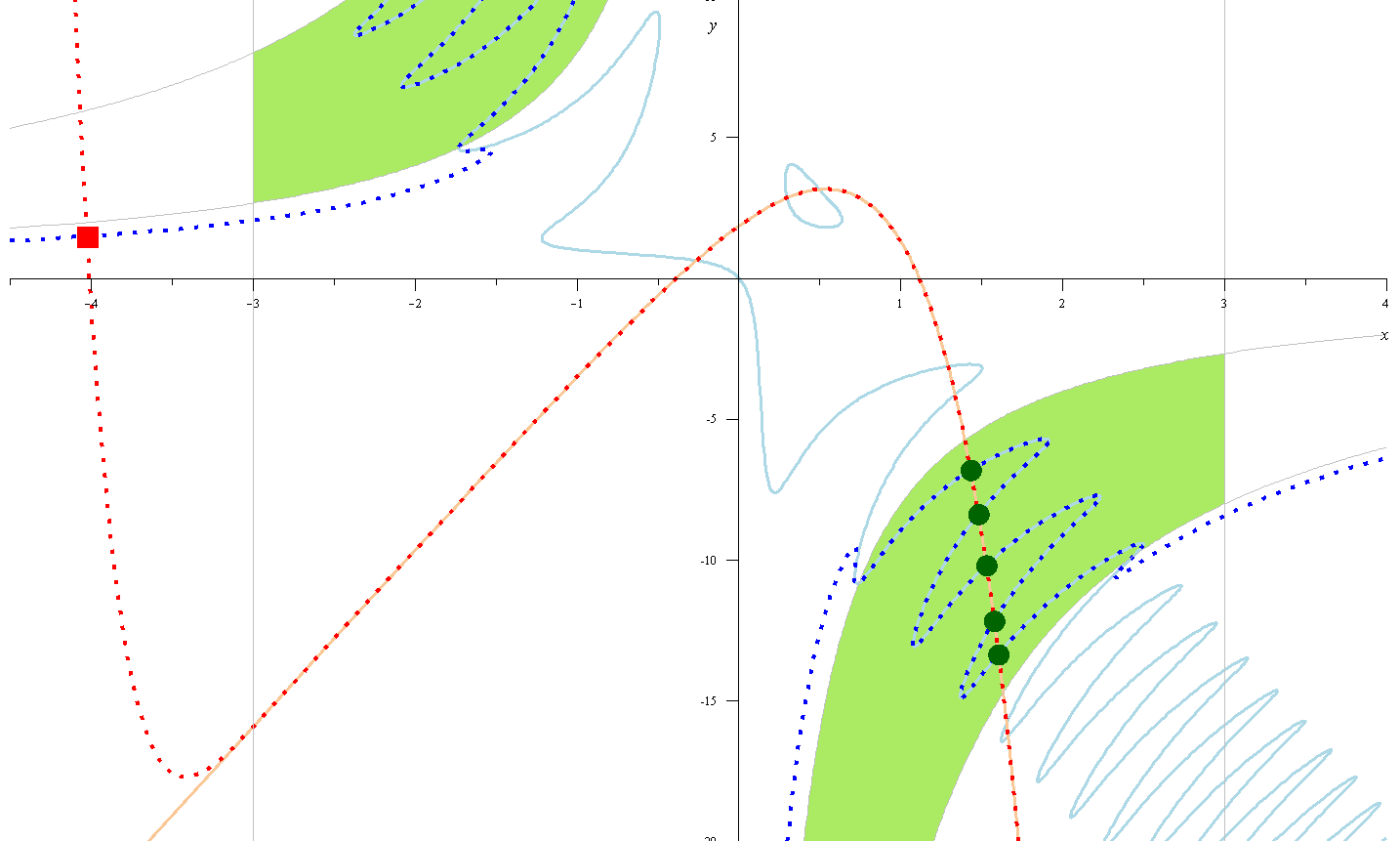}
\end{center}
\caption{Solutions to the polynomial system $p_2(x,y) = 0, q_2(x,y) = 0$}
\end{figure}

Continuing the procedure, we will find another two roots in search region 3.

\end{example}

\begin{example}
Let us consider the following system of nonlinear equations studied by Krzyworzcka \cite[Example 2]{Krzyworzcka1996}:
\begin{eqnarray*}
-x_1 - 0.75 - 0.25 x_2^2 x_4 x_6 = 0, \\
-x_2 - 0.405 e^{1+x_1 x_2} + 1.405 = 0, \\
-x_3 + 0.5 x_4 x_6 -1.5 = 0, \\
-x_4 + 0.605 e^{1-x_3^2} + 0.395 = 0, \\
-x_5 + 0.5 x_2 x_6 - 1.5 = 0, \\
-x_6 + x_1 x_5 = 0. \\
\end{eqnarray*}

Suppose that no additional information is given on the estimated location of the root(s).
Write the Chebysev approximation of the transcendental function $e^y$

\begin{center}
$
 3.4209729272\cdot10^{-12}y^{15}
+5.6225114884\cdot10^{-11}y^{14}
-4.1022516049\cdot10^{-10}y^{13}
-6.9288841558\cdot10^{-9}y^{12}
+8.3636655828\cdot10^{-8}y^{11}
+0.1122572549\cdot10^{-5}y^{10}
-5.3510542091\cdot10^{-7}y^9
-0.0000174769 y^8
+0.0003012363 y^7
+0.0025115198 y^6
+0.0066485671 y^5
+0.0271995805 y^4
+0.1789743081 y^3
+0.5704984723 y^2
+0.9738974605 y
+0.9439171466
$
\end{center}

that approximates the exponential function with accuracy 0.01 in the interval $-10 \leq y \leq 10$.
First we search for the solutions satisfying $ -10 \leq 1+x_1 x_2 \leq 10$ (from the second equation) and
$ -10 \leq 1-x_3^2 \leq 10$ (from the fourth equation). The approximating polynomial system has 1343
solutions, found by HOM4PS-3, four of them real.

The first solution of the polynomial system is
\begin{eqnarray*}
x_1 = -0.980743705883039493312476489776423 \\
x_2 = 1.02465410949677159203773809309151 \\
x_3 = -1.06045303397770105733015703750029 \\
x_4 = .897825333984771728771461332720686 \\
x_5 = -0.998361665678699229116492980275557 \\
x_6 = .979136919809291600947871819079083
\end{eqnarray*}
and the Newton-iteration written for the original system, starting from this point, converges to the near point
\begin{eqnarray*}
x_1 = -1 \\
x_2 =  1  \\
x_3 = -1  \\
x_4 =  1 \\
x_5 = -1  \\
x_6 =  1
\end{eqnarray*}
that was also found by Krzyworzcka \cite[Example 2]{Krzyworzcka1996}.

The second solution of the polynomial system is
\begin{eqnarray*}
x_1 = -1.03416093043463753827399238895323 \\
x_2 =  -0.545512751309159592218439384752022 \\
x_3 =  0.409785310325564054403836837803004 \\
x_4 =  1.76772558666096248026600321202482 \\
x_5 =  -2.08935179017993586914288645096315 \\
x_6 =  2.16072599133775806444533214652971
\end{eqnarray*}
and the Newton-iteration written for the original system, starting from this point, converges to the near point
\begin{eqnarray*}
x_1 = -1.04320094527705867310938476801 \\
x_2 = -.550936201394894459163931845959 \\
x_3 =  0.431936026252448558564000278415 \\
x_4 =  1.75965881945053520580861011514 \\
x_5 = -2.10487492454578699747973642815  \\
x_6 = 2.19580751097614254519544578876
\end{eqnarray*}
that was also found by Ibrahim and Tawhid \cite[Table 2, column PSO]{IbrahimTawhid2019}.

It worths noting that although the third and the fourth solutions of the polynomial system
\begin{eqnarray*}
x_1 = -0.444835925873194521485096536120813 \\
x_2 = 0.562040139916147599444071427709663 \\
x_3 = -3.43209614317652617330132315819424 \\
x_4 = -6.51513108377923871925967726518616 \\
x_5 = -1.33332406168373273070171870543547 \\
x_6 = 0.593110443468091572740967278548609
\end{eqnarray*}
and
\begin{eqnarray*}
x_1 = -2.68376423295794954961057650788998 \\
x_2 = 4.34909015091502538174568080077468  \\
x_3 = -1.29552665643317582956675523161942 \\
x_4 =  0.694434018201516277085069445214144 \\
x_5 = -0.219427647087008066805129831853443 \\
x_6 = 0.588892070974231857286181864692679
\end{eqnarray*}
do not satisfy our assumptions $ -10 \leq 1+x_1 x_2 \leq 10$ (from the second equation) or
$ -10 \leq 1-x_3^2 \leq 10$ (from the fourth equation), the Newton-iteration, starting from these points,
still converges to $[-1,1,-1,1,-1,1]$ again, which is not at all close to any of the two points.

\end{example}

\section*{}
Since the aim of this paper is only to illustrate, by small examples, that paxotopy method is capable of solving
nonlinear systems, a detailed comparison to other methods, such as the
multiobjective evolutionary algorithm (see e.g. \cite{GrosanAbraham2008}),
adaptive simulated annealing \cite{OliveiraIngberPetragliaPetragliaMachado2012},
hybrid methods of cuckoo search and particle swarm optimization \cite{IbrahimTawhid2019},
Lanczos and conjugate gradient squared methods \cite{Krzyworzcka1996},
chaos optimization and quasi-Newton method \cite{LuoTangZhou2008},
continuous greedy randomized adaptive search procedure \cite{HirschPardalosResende2009},
biased random-key genetic algorithm \cite{SilvaResendePardalos2014}
is a subject of future research.

\bibliographystyle{plainnat}
\bibliography{paxotopy-references}

@book{Boyd2014,
  author    = {Boyd, J. P.},
  title     = {Solving Transcendental Equations: The Chebyshev Polynomial Proxy and Other Numerical Rootfinders, Perturbation Series, and Oracles},
  publisher = {SIAM},
  year      = {2014}
}

@article{Bozoki2020,
  author  = {Boz'{o}ki, S.},
  title   = {Eccentric pie charts and an unusual pie cutting},
  journal = {Information Visualization},
  volume  = {19},
  number  = {4},
  pages   = {288--295},
  year    = {2020},
  doi     = {10.1177/1473871620925078}
}

@book{Carnahan1964,
  author    = {Carnahan, B.},
  title     = {Applied Numerical Methods (Preliminary edition)},
  publisher = {Wiley},
  year      = {1964}
}

@book{CarnahanGourdinLutherWilkes1969,
  author    = {Carnahan, B. and Gourdin, A. and Luther, H. A. and Wilkes, J. O.},
  title     = {Applied Numerical Methods},
  publisher = {Wiley},
  year      = {1969}
}

@article{GrosanAbraham2008,
  author  = {Grosan, C. and Abraham, A.},
  title   = {Multiple solutions for a system of nonlinear equations},
  journal = {International Journal of Innovative Computing, Information and Control},
  volume  = {4},
  number  = {9},
  pages   = {2161--2170},
  year    = {2008}
}

@article{HirschPardalosResende2009,
  author  = {Hirsch, M. J. and Pardalos, P. M. and Resende, M. G. C.},
  title   = {Nonlinear Analysis: Real World Applications},
  journal = {Nonlinear Analysis: Real World Applications},
  volume  = {10},
  number  = {4},
  pages   = {2000--2006},
  year    = {2009},
  doi     = {10.1016/j.nonrwa.2008.03.006}
}

@article{IbrahimTawhid2019,
  author  = {Ibrahim, A. M. and Tawhid, M. A.},
  title   = {A hybridization of cuckoo search and particle swarm optimization for solving nonlinear systems},
  journal = {Evolutionary Intelligence},
  year    = {2019},
  doi     = {10.1007/s12065-019-00255-0}
}

@article{JiWuFengLiQin2016,
  author  = {Ji, Z. and Wu, W. and Feng, Y. and Li, Y. and Qin, X. L.},
  title   = {Numerical Method for Real Root Isolation of Semi-Algebraic System and Its Applications},
  journal = {Journal of Computational and Theoretical Nanoscience},
  volume  = {13},
  number  = {1},
  pages   = {803--811},
  year    = {2016},
  doi     = {10.1166/jctn.2016.4878}
}

@article{Krzyworzcka1996,
  author  = {Krzyworzcka, S.},
  title   = {Extension of the Lanczos and CGS methods to systems of nonlinear equations},
  journal = {Journal of Computational and Applied Mathematics},
  volume  = {69},
  number  = {1},
  pages   = {181--190},
  year    = {1996},
  doi     = {10.1016/0377-0427(95)00032-1}
}

@article{LuoTangZhou2008,
  author  = {Luo, Y.-Z. and Tang, G.-J. and Zhou, L.-N.},
  title   = {Hybrid approach for solving systems of nonlinear equations using chaos optimization and quasi-Newton method},
  journal = {Applied Soft Computing},
  volume  = {8},
  number  = {2},
  pages   = {1068--1073},
  year    = {2008},
  doi     = {10.1016/j.asoc.2007.05.013}
}

@incollection{OliveiraIngberPetragliaPetragliaMachado2012,
  author    = {Oliveira Jr., H. and Ingber, L. and Petraglia, A. and Petraglia, M. R. and Machado, M. A. S.},
  title     = {Nonlinear Equation Solving},
  booktitle = {Stochastic Global Optimization and Its Applications with Fuzzy Adaptive Simulated Annealing},
  series    = {Intelligent Systems Reference Library},
  volume    = {35},
  chapter   = {10},
  pages     = {169--187},
  publisher = {Springer},
  address   = {Berlin, Heidelberg},
  year      = {2012},
  doi       = {10.1007/978-3-642-27479-4_10}
}

@article{SilvaResendePardalos2014,
  author  = {Silva, R. M. A. and Resende, M. G. C. and Pardalos, P. M.},
  title   = {Finding multiple roots of a box-constrained system of nonlinear equations with a biased random-key genetic algorithm},
  journal = {Journal of Global Optimization},
  volume  = {60},
  number  = {2},
  pages   = {289--306},
  year    = {2014},
  doi     = {10.1007/s10898-013-0105-7}
}


\end{document}